\documentclass{amsart}

\usepackage{overpic}
\usepackage{graphicx}%
\usepackage{multirow}%
\usepackage{amsmath,amssymb,amsfonts}%
\usepackage{amsthm}%
\usepackage{mathrsfs}%
\usepackage[title]{appendix}%
\usepackage{xcolor}%
\usepackage{textcomp}%
\usepackage{manyfoot}%
\usepackage{booktabs}%
\usepackage{algorithm}%
\usepackage{algorithmicx}%
\usepackage{algpseudocode}%
\usepackage{listings}%
\usepackage{mathrsfs}

\usepackage{tikz}
\usepackage{tikz-3dplot}
\usetikzlibrary{knots}

\usepackage{amsmath}
\usepackage{amsfonts}
\usepackage{amssymb}
\usepackage{graphicx} 
\usepackage{overpic}

\usepackage{tikz}
\usepackage{tikz-3dplot}
\usetikzlibrary{knots}

\theoremstyle{plain}
\newtheorem{thm}{Theorem}[section]

\newtheorem{lemma}[thm]{Lemma}
\newtheorem{prop}[thm]{Proposition}
\newtheorem{claim}[thm]{Claim}

\theoremstyle{remark}
\newtheorem{remark}[thm]{Remark}

\begin{document}

\title[Growth function of the cactus groups]{Growth function of the cactus groups}

\author{Takatoshi Hama}
\address{Graduate School of Integrated Basic Sciences, Nihon University, 3-25-40 Sakurajosui, Setagaya-ku, Tokyo 156-8550, Japan}
\email{chta23018@g.nihon-u.ac.jp}


\date{\today}



\begin{abstract}
We determine the growth functions of the cactus groups of all degrees with respect to their standard generators, and provide explicit equations for specific degrees as examples. 
\end{abstract}

\maketitle 

\section{Introduction}
The cactus group $J_n$ was originally introduced in the context of a real form of the Deligne-Mumford space $\overline{M}_{0,n+1}$ of stable genus $0$ curves with $n+1$ marked points \cite{devadoss2000tessellationsmodulispacesmosaic}. 
On the other hand, research on this group from the view point of geometric group theory has been conducted by Genevois in \cite{genevois2022cactusgroupsviewpointgeometric}. 
Furthermore, the cactus group can be interpreted as a group possessing characteristics of both the coxeter group and the braid group. In fact, it is known that the cactus group acts on geometric objects. See \cite{devadoss2000tessellationsmodulispacesmosaic}, \cite{genevois2022cactusgroupsviewpointgeometric}. 
Additionally, in geometric group theory, the calculation of the growth function of the coxeter group has been extensively studied. Scott gives an extension of the growth function of the coxeter group in \cite{Scott2007GrowthSF}. The following is obtained by using these results. 

\begin{thm}\label{growth_of_cactus_gr}
    Let $n$ be a positive integer greater than $2$ and $\{a_i\}_{i \in \mathbb{N}}$ be the sequence satisfying the following equation:
    \[ \displaystyle
    \left( \sum_{i=0}^{\infty} a_i x^i \right) \left( \sum_{k=0}^{n-1} d_{n,k} (-x)^k (1+x)^{n-k-1}  \right) = (1+x)^{n-1}, 
    \]
    where 
\begin{align*} \displaystyle
 d_{n,k} := D_{n+1,k} + D_{n+1,k-1}
\end{align*}
    Here, $D_{n+1,k-1}$ is the Kirkman-Cayley number, that is, equal to 
    \[
    \displaystyle
 D_{n+1,k-1} =\frac{1}{k} \begin{pmatrix}
 n-2 \\ k-1 
\end{pmatrix}
\begin{pmatrix}
 n+k-1 \\ k -1 
\end{pmatrix}. 
    \]
    Then, the growth function $g_{(J_n, S)}$ of the cactus group $J_n$ of degree $n$ with respect to the standard generator S is given as follows. 
    \[
    g_{(J_n, S)}(m) = \sum_{i=0}^{m} a_i 
    \]
\end{thm}

\section{Cactus groups and geometric group theory}
In this section, we give the purely algebraic definition of the cactus group, and some basic definitions and propositions of the geometric group theory. We follow the notation for the cactus group used in \cite{genevois2022cactusgroupsviewpointgeometric}. 

For any integer $n \geq 2$, the \textit{cactus group of degree $n$}, denoted by $J_n$, is given by the presentation with the generators $s_{p,q}$ for $1 \leq p < q \leq n$, subject to the following relations. 
\begin{itemize}
   \item $s_{p,q}^2 = e$ for every $1 \leq p < q \leq n$, 
   \item $s_{p,q}s_{m,r} = s_{m,r}s_{p,q}$ for all $1 \leq p < q \leq n$ and $1 \leq m < r \leq n$ satisfying $[p, q] \cap [m, r] = \emptyset$, 
   \item $s_{p,q}s_{m,r} = s_{p+q-r,p+q-m}s_{p,q}$ for all $1 \leq p < q \leq n$ and $1 \leq m < r \leq n$ satisfying $[m, r] \subset [p, q]$. 
\end{itemize}
Here, $e$ denotes the identity element, and $[p, q]$ denotes the discrete interval ${ p, p+1, \dots, q-1, q }$ of integers for positive integers $p$ and $q$ with $p < q$. 

We also introduce a subgroup that will be used in this paper. 
For every integer $n \geq 2$ and subset $S \subset [2, n]$, let $J_n^S$ be the subgroup of $J_n$ generated by the elements $s_{p,q}$ for $1 \leq p < q \leq n$ with $q - p + 1 \in S$, and defined by the following relations. 
\begin{itemize}
   \item $s_{p,q}^2 = e$ for every $1 \leq p < q \leq n$ satisfying $q-p+1\in S$, 
   
   \item $s_{p,q}s_{m,r} = s_{m,r}s_{p,q}$ for every $1 \leq p < q \leq n$ and $1 \leq m < r \leq n$ satisfying $q-p+1\in S$, $r-m+1 \in S$, and $[p, q] \cap [m, r] = \emptyset$. 
   
   \item $s_{p,q}s_{m,r} = s_{p+q-r,p+q-m}s_{p,q}$ for every $1 \leq p < q \leq n$ and $1 \leq m < r \leq n$ satisfying $q-p+1\in S$, $r-m+1 \in S$, and $[m, r] \subset [p, q]$. 
\end{itemize}

In the subsequent sections, our primary focus is the subgroup $J_n^{[2,n-1]} \subset J_n$, generated by the generators of $J_n$ excluding $s_{1,n}$. The defining relations are given by those of $J_n$ that do not contain $s_{1,n}$. 
It is known that there exists a semi-direct product decomposition for $J_n$:
\begin{prop}\label{semidirectprod}
    The cactus group $J_n$ admits a semidirect product decomposition:$J_n \cong J_n^{[2,n-1]} \rtimes J_n^{\{n\}}$ where $J_n^{\{n\}} = \langle s_{1,n} \mid s_{1,n}^2 = 1 \rangle \cong \mathbb{Z}/2\mathbb{Z}$. 
\end{prop}

\begin{proof}
By \cite[Corollary~5.5]{genevois2022cactusgroupsviewpointgeometric}, we can obtain this proposition directly.     
\end{proof}

We next recall some basic notions of geometric group theory. 

A \textit{cubical complex} or \textit{cube complex} $X$ is the quotient of a disjoint union of cubes 
\[ \displaystyle
K = \bigsqcup_{\lambda
 \in \Lambda } I^{n_\lambda}
\]
by an equivalent relation such that 
the restrictions $p_\lambda: I^{n_\lambda} \longrightarrow X$ of the natural projection $p: K \longrightarrow X$ are required to satisfy: 
\begin{itemize}
    \item $p_\lambda$ is injective for any $\lambda$ in $\Lambda$. 
    
    \item If $p_\lambda (I^{n_\lambda}) \cap p_{\lambda'} (I^{n_{\lambda'}}) \neq \emptyset$, then there is isometry $h_{\lambda, \lambda'}$ from a face $T_{\lambda}$ contained in $I^{n_\lambda}$ onto a face $T_{\lambda'}$ contained in $I^{n_{\lambda'}}$ such that $p_\lambda(t_{ \lambda }) = p_{\lambda '}(h_{\lambda, \lambda'}(t_{ \lambda }))$ for all $t_{\lambda}$ in $T_{\lambda}$. 
\end{itemize}

For an $k$-cube $I^k$ and its vertex $x$ in $(I^k)^{(0)}$, the \textit{geometric link} Lk$(x,I^k)$ of $x$ in $I^k$ is the set of unit vectors at $x$ that point into $I^k$. 
For a cubical complex $X$ and its vertex $x$ in $X^{(0)}$, the \textit{geometric link} of $x$ in $X$, denoted by Lk$(x,X)$, is the space of unit initial vectors at $x$ endowed with the quotient metric associated to the projection 
\[ \displaystyle
\bigsqcup_{\lambda \in \Lambda} \text{Lk}(x_{\lambda},(I^{k})_{\lambda}) \longrightarrow \text{Lk}(x,X) 
\]
induced by 
\[ \displaystyle
\bigsqcup_{\lambda \in \Lambda} (I^{k})_{\lambda} \longrightarrow X 
\]
\begin{remark}
    We remark that the geometric link of a vertex in a cubical complex is a simplicial complex. See \cite{BridsonHaefliger} for more details. 
\end{remark}

From the viewpoint of geometric group theory, the following result was established about the cactus group in \cite{genevois2022cactusgroupsviewpointgeometric}. 
Note that we denote Cay$(G,S)$ by the Cayley graph of a group $G$ with respect to generating set $S$ in this paper. 
\begin{prop}[{\cite[Theorem~1.1]{genevois2022cactusgroupsviewpointgeometric}}] \label{ccts median}
For every $n \geq 2$, the Cayley graph Cay$(J_n,S)$ (denoted by $\mathscr{C}_n$ in \cite{genevois2022cactusgroupsviewpointgeometric}) of $J_n$ with respect to the standard generating set $S$ is median. 
\end{prop}

Consequently, the cube-completion of Cay$(J_n,S)$ forms a CAT$(0)$ cubical complex by the next proposition, where CAT$(0)$ space is defined as follows. 

A metric space $(X, d)$ is called a \textit{CAT(0) space} if every geodesic triangle $(\gamma_0: [0,L_1] \to X, \gamma_1: [0,L_2] \to X, \gamma_2: [0,L_3] \to X)$ in $X$ satisfies the following condition: 
For a comparison triangle $(\gamma_0' : [0,L_1] \to \mathbb{R}^2, \gamma_1' : [0,L_1] \to \mathbb{R}^2, \gamma_2' : [0,L_2] \to \mathbb{R}^2)$ in the Euclidean plane $(\mathbb{R}^2, d_{\mathbb{R}^2})$ with the same side lengths, for any indices $j, k \in \{0, 1, 2\}$, and for any $s \in [0, L_j]$ and $t \in [0, L_k]$, the following inequality holds. 
\[
d(\gamma_j (s), \gamma_k (t)) \leq d_{\mathbb{R}^2} (\gamma_j' (s), \gamma_k' (t)) 
\]

\begin{prop}[{\cite[Theorem~2.3]{genevois2022cactusgroupsviewpointgeometric}}]\label{Median_leads_CAT0}
 A graph is median if and only if its cube completion is $CAT(0)$. 
\end{prop}
Here, the \textit{cube-completion of a median graph} is the cubical complex obtained by filling in all the one-skeleton of cubes in the graph with cubes. 
See \cite{genevois2022cactusgroupsviewpointgeometric} and \cite[Theorem~6.1]{ChepoiCube} for details. 

Furthermore, the following well-known proposition deduces that $J_n$ acts the cube-completion of $\mathscr{C}_n$ cocompactly by the left $G$-action. 
\begin{prop} \label{cocompact action}
    Let $G$ be a group generated by a finite generating set $S$ and Cay$(G,S)$ be the Cayley graph of $G$ with respect to $S$. 
    Then the left $G$-action on Cay$(G,S)$ is a cocompact action. 
\end{prop}
For completeness, we give an elementary proof. 
    \begin{proof}
    The action is transitive on the vertices, since every vertex $g$ equals $g\cdot e$. 
    Every edge of Cay$(G,S)$ is of the form $\{g,gs\}=g\cdot\{e,s\}$ for some 
    $g\in G$ and $s\in S$, so there are at most $|S|$ orbits of edges. 
    Since the quotient graph has one vertex and at most $|S|$ edges, the quotient graph is 
    a bouquet consisting of at most $|S|$ circles. It is a finite graph, 
    so the left $G$-action is cocompact. 
\end{proof}

For a CAT$(0)$ cubical complex $X$ with a cocompact cellular action by a group $G$, the next proposition is obtained in \cite{Scott2007GrowthSF}. 

\begin{prop}[{\cite[Theorem~1]{Scott2007GrowthSF}}]\label{growth_series}
Let $x$ be a vertex in a cubical complex $X$. 
If the links of all vertices have the same $f$-polynomial $f_x$ of $x$, then 
\[
G_x(t) f_x \left( \frac{-t}{1+t} \right) = 1. 
\]
\end{prop}

Here, we recall the definitions of the $f$-polynomial, the full growth series, and the growth function. See \cite{Scott2007GrowthSF} for more details. 

The \textit{$f$-polynomial} of a simplicial complex $L$ is given by
\[
f_L(t) = \sum_{\sigma \in L} t^{\dim \sigma + 1}.
\]
Note that we assume that $L$ contains an empty simplex of dimension $-1$, so the $f$-polynomial always has the constant term $1$. For a vertex $x$ in $X^{(0)}$, $f_x$ denotes the $f$-polynomial of the  geometric link of the vertex $x$ in $X^{(0)}$. 

For a CAT$(0)$ cubical complex with a cellular action by a group $G$, the \textit{full growth series} $G_x(t)$ of $X$ as seen from a vertex $x$ in $X$ is defined by the following: 
\[
G_x (t)= \sum_{y \in X} t^{d(x,y)}, 
\]
where $d(x, y)$ is the $l_1$-distance between vertices $x$ and $y$ of $X$. 

For a group $G$ generated by a finite generating set $S$, the \textit{growth function} $f_{(G,S)}(n)$ of $G$ with respect to $S$ is the map from $\mathbb{N}$ to $\mathbb{N}$ given as follows: 
\[ \displaystyle 
g_{(G,S)}(n) := \sum_{i=0}^{n} c_i, 
\]
where $c_i$ is the number of the elements whose word length is equal to $i$. 

We note that if the cube-completion $C$ of the Cayley graph of $G$ with respect to $S$ is a CAT$(0)$ cubical complex, the sum of coefficient of the full growth series of $C$ gives the growth function of $G$ with respect to $S$. 

\section{Parenthesization and the Kirkman-Cayley number}
In this section, we recall the definition of the parenthesization and a well known proposition about it. See \cite{GAIFFI2015279} for more details. 

Let $a_1$, $a_2$, $\cdots$, $a_n$ be a sequence (called ``list" in \cite{GAIFFI2015279}) of distinct numbers. A \textit{parenthesization} of the sequence $a_1$, $a_2$, $\cdots$, $a_n$ is a set $\mathscr{P}$ of discrete intervals $[i,j] = \{ i,\ i+1,\ \cdots j \}$ in $\{ 1,\ 2,\ \cdots,\ n \}$ with $1 \leq i < j \leq n$ satisfying the following conditions. 
\begin{itemize}
    \item $\mathscr{P}$ contains $[1,n]$. 
    \item For any intervals $I$, $J$ in $\mathscr{P}$, either $I \subset J$, $J \subset I$, or $I \cap J = \emptyset$. 
\end{itemize}
Each interval in $\mathscr{P}$ is called a \textit{couple of parentheses}. We denote by $\mathscr{P}_2(( a_1, a_2, \cdots, a_n),k)$ the set of parenthesizations with $|\mathscr{P}|=k$. 
We remark that we use the different notation from \cite{GAIFFI2015279} to avoid confusing the notation for a generating set. 

We next recall the Kirkman-Cayley number $D_{n+1,k-1}$. 

Let $P_{n+1}$ be a convex polygon on the Euclidean plane $\mathbb{R}^2$ with $n+1$ vertices labeled counterclockwise as $v_0, v_1, \dots, v_n$. 
For each $i \in \{1, 2, \dots, n\}$, let $e_i$ denote the boundary edge $(v_{i-1}, v_i)$, and let $e_0$ denote the base edge $(v_n, v_0)$. 

We denote by $\mathcal{D}_{n+1, k-1}$ the set of dissections of $P_{n+1}$ formed by $k-1$ mutually non-intersecting diagonals in its interior. 

Set $D_{n+1, k-1} := |\mathcal{D}|_{n+1, k-1}$. 
Here, $D_{n+1, k-1}$ is called the \textit{Kirkman-Cayley number}. 
We remark that if $n-k \leq 1$, $D_{n+1,k-1}$ is set to be $0$. See \cite{przytycki2000polygon} for more details. 

As is well known, there exists a bijection between $\mathscr{P}_2(( a_1, a_2, \cdots, a_n),k)$ and $\mathcal{D}_{n+1, k-1}$. 

\begin{lemma}
    $|\mathscr{P}_2(( 1, 2, \cdots, n),k)|$ coincides with the Kirkman-Cayley number $D_{n+1,k-1}$. 
\end{lemma}

\begin{proof} \label{K-C coincides with parenthesization number}
We first establish a local rule that maps a couple of parentheses to a unique diagonal of the polygon. 

Let $p = [i, j]$ be a couple of parentheses in $\mathscr{P}$. 
If $p \neq [1,n]$, we define the directed diagonal $f(p)$ as the line segment connecting the starting vertex $v_{i-1}$ of edge $e_i$ and the ending vertex $v_j$ of edge $e_j$: 
\[
f(p) := (v_{i-1}, v_j). 
\]
Having established the local assignment $f$, we naturally extend it to the entire parenthesization structure.

Let $\mathscr{P}_2(( 1, 2, \cdots, n),k)$ be a set of parenthesizations. 
The map 
\[
\Psi:\mathscr{P}_2(( 1, 2, \cdots, n),k)  \longrightarrow \mathcal{D}_{n+1, k-1}
\]
is defined by applying $f$ to each non-maximal parenthesis $p$ in the parenthesization $P$ in 
$\mathscr{P}_2(( 1, 2, \cdots, n),k)$: 
\[
\Psi(P) := \left\{ f(p) \mid p \in P \setminus \{[1,n]\}, P  \in \mathscr{P}_2(( 1, 2, \cdots, n),k) \right\}. 
\]
Therefore, $|\mathscr{P}_2(( 1, 2, \cdots, n),k)|$ coincides with the Kirkman-Cayley number $D_{n+1,k-1}$. 
\end{proof}

Furthermore, an explicit equation is given for the Kirkman–Cayley numbers. 

\begin{prop}[{\cite[Corollary~4.1]{GAIFFI2015279}}]\label{K_C_number}
    \[
    \displaystyle
 D_{n+1,k-1} =\frac{1}{k} \begin{pmatrix}
 n-2 \\ k-1 
\end{pmatrix}
\begin{pmatrix}
 n+k-1 \\ k -1 
\end{pmatrix}. 
    \]
\end{prop}

\section{Proof of Theorem~1.1}

In this section, we give a proof of Theorem~\ref{growth_of_cactus_gr}. 

\begin{proof}[Proof of Theorem~\ref{growth_of_cactus_gr}]
    Let $J_n$ be the cactus group of degree $n$ with respect to the standard generating set $S$. 

    By Proposition~\ref{ccts median}, Proposition~\ref{Median_leads_CAT0} and Proposition~\ref{cocompact action}, we conclude that $J_n$ acts on the cube-completion $\mathscr{A}_n$ of Cay$(J_n,S)$ as cocompact cellular action. 
    Furthermore, the action on $\mathscr{A}_n$ by $J_n$ is transitive. It concludes the geodesic links of all vertices of $\mathscr{A}_n$ have the same $f$-polynomial. 
    Therefore, we can apply Proposition~\ref{growth_series} for $J_n$ and $\mathscr{A}_n$. 

    We next determine the $f$-polynomial $f_e(x)$ of the geodesic link around $e$ in $\mathscr{A}_n$ with cellular action by $J_n$.
    Thus, we count the numbers of the simplices with same dimension in the geometric link around $e$. 
    Here, recall that the numbers of the $k$-cubes around $e$ in $\mathscr{A}_n$ coincides with the numbers of $k-1$-simplices in the geodesic link of $e$ in $\mathscr{A}_n$. 
   
    Let $\mathcal{C}^{[2,n-1]}_k$ the set of $k$-cubes in $\mathscr{A}^{[2,n-1]}$ around the vertex $e$, where $\mathscr{A}^{[2,n-1]}$ is the cube-completion of the Cayley graph  $\operatorname{Cay}(J_n^{[2,n-1]}, S \setminus \{s_{1,n}\})$ of $J_n^{[2,n-1]}$. 
    That is, $\mathcal{C}^{[2,n-1]}_k$ can be regarded as the set of the subset 
    $ \{ s_{i(1),j(1)}, \cdots, s_{i(k),j(k)} \}$ of the generating set $S \setminus \{ s_{1,n} \}$ satisfying the next conditions: for all $l,m$ in $\{ 1,\cdots , k \}$, 
    $[i(l),j(l)] \subset [i(m),j(m)]$ or $[i(l),j(l)] \supset [i(m),j(m)]$ or $[i(l),j(l)] \cap [i(m),j(m)] = \emptyset$. 
    Let $\operatorname{Index}(\mathcal{C}^{[2,n-1]}_k)$ be the set of the set of intervals whose endpoints are the index of the generators in $\mathcal{C}^{[2,n-1]}_k$: 
    \[
    \operatorname{Index}(\mathcal{C}^{[2,n-1]}_k) :=
    \{
    \{ [i(1),j(1)], \cdots , [i(k),j(k)] \} \mid
    \{ s_{i(1),j(1)}, \cdots , s_{i(k),j(k)} \} \in \mathcal{C}^{[2,n-1]}_k
    \}
    \]
    We can also describe this set as follows. 
    \[
    \operatorname{Index}(\mathcal{C}^{[2,n-1]}_k) = 
    \{
    \{ I_1, \cdots I_k \} \mid  \forall \alpha , \beta \in \{ 1, \cdots k \}, I_\alpha \subset I_\beta \text{ or }  I_\alpha \supset I_\beta \text{ or } I_\alpha \cap I_\beta = \emptyset 
    \}
    \]
    Clearly, there exists a bijection between $\mathcal{C}^{[2,n-1]}_k$ and $\operatorname{Index}(\mathcal{C}^{[2,n-1]}_k)$. 

    On the other hand, 
    $\mathscr{P}_2(( 1, 2, \cdots, n),k)$ is the set $\{ I_1, \cdots I_k \}$ of the subintervals in $[1,n]$ satisfying the following conditions. 
    \begin{itemize}
        \item $[1,n] \in \{ I_1, \cdots , I_k \}$. 
        \item For all $\alpha$ and $\beta$, $I_\alpha \subset I_\beta$ or $I_\alpha \supset I_\beta$ or$ I_\alpha \cap I_\beta = \emptyset$. 
    \end{itemize}
    We can describe $\mathscr{P}_2(( 1, 2, \cdots, n),k)$ by the following set. 
    \[
    \{
    \{ I_1, \cdots , I_{k-1} \} \sqcup \{ [1,n] \} \mid \forall \alpha , \beta \in \{ 1, \cdots k \}, I_\alpha \subset I_\beta \text{ or }  I_\alpha \supset I_\beta \text{ or } I_\alpha \cap I_\beta = \emptyset
    \}
    \]
    Therefore, the number of $\mathcal{C}^{[2,n-1]}_k$ and the number of $\mathscr{P}_2(( 1, 2, \cdots, n),k)$ coincide. 
    By Lemma~\ref{K-C coincides with parenthesization number}, the number of $\mathcal{C}^{[2,n-1]}_k$ coincides with $D_{n+1,k}$. 

    By Proposition~\ref{semidirectprod}, it is implied that Cay$(J_n, S)$ is isometric to $\operatorname{Cay}(J_n^{[2,n]}, S \setminus \{s_{1,n}\}) \times \{0,1\}$ as a graph, where the two copies are connected by edges corresponding to $s_{1,n}$. 
    Therefore, the geodesic link of $e$ in $\mathscr{A}_n$ is the cone of the geodesic link of $e$ in $\mathscr{A}_n^{[2,n-1]}$. 
    Since the cube around $e$ in $\mathscr{A}_n^{[2,n-1]}$ does not have $s_{1,n}$ as an edge, the number of the $k$-cube around $e$ in $\mathscr{A}_n$ is equal to the sum of the number of $\mathcal{C}^{[2,n-1]}_k$ and $\mathcal{C}^{[2,n-1]}_{k-1}$. 
    Therefore, the number of cubes around $e$ in $\mathscr{A}_n$ coincides with $D_{n+1,k} + D_{n+1,k-1}$. 
    It concludes that the number of the $k-1$-simplices in the geodesic link of $e$ in $\mathscr{A}_n$ coincides with $D_{n+1,k} + D_{n+1,k-1}$. 
    Furthermore, the following claim concludes that the highest dimension of the cubes in $\mathscr{A_n}^{[2,n-1]}$ is $n-2$. 
    \begin{claim}[{\cite[Proposition~7.2]{genevois2022cactusgroupsviewpointgeometric}}] \label{cube dimension}
    Fix an $n \ge 3$. The cube-completion $\mathscr{A}_n^{[2,n-1]}$ of the median graph $\mathscr{C}_n^{[2,n-1]}$ is an $(n - 2)$-manifold.
    \end{claim}
    \noindent
    It implies the highest dimension of the cubes in $\mathscr{A_n}$ is $n-1$. 
    
    Therefore the $f$-polynomial $f_e (t)$ of Lk$(e,\mathscr{A}_n)$ is given by 
    \[ \displaystyle
    f_e (t) = \sum_{k=0}^{n-2} d_{n,k} t^k, 
    \]
    where $d_{n,k} := D_{n+1,k} + D_{n+1,k-1}$. 

    By Proposition~\ref{growth_series}, the following holds. 
    \[ \displaystyle
    G_e (t) \sum_{k=0}^{n-1} d_{n,k} \left(\frac{-x}{1+x} \right)^k = 1
    \]
    Simplifying the above equation, we can obtain 
    \[\displaystyle
    G_e (t) \left( \sum_{k=0}^{n-1} d_{n,k} (-x)^k (1+x)^{n-k-1}  \right) = (1+x)^{n-1}
    \]
    Let $\{ a_i \}_{i \in \mathbb{N}}$ be the sequence of coefficient of $G_e (t)$. Then $\{ a_i \}_{i \in \mathbb{N}}$ gives the growth function $g_{(J_n,S)}(m)$ by
    \[
    g_{(J_n, S)}(m) = \sum_{i=0}^{m} a_i.  
    \]

\end{proof}

\section{Examples}
In this section, we give the full growth series of the cube-completion of $\mathscr{A}_n$ for specific cases of degree $3$, $4$, $5$, and $6$. 

\bigskip

When $n=3$, the following holds. 

\begin{align*}
    d_{3,0} = 1 \\
    d_{3,1} = 3\\ 
    d_{3,2} = 2
\end{align*}

Therefore, the full growth series $G_e(t)$ of the cube-completion of $\mathscr{A}_3$ satisfies the following. 
\begin{align*} \displaystyle
    G_e (t) &= \frac{(1+t)^2}{ (1+t)^2 - 3t(1+t)^1 + 2t^2(1+t)^0} \\
&= \frac{(1+t)^2}{1-t}
\end{align*} 
Then we obtain the next recurrence relation of the sequence of the coefficients by applying the power series expansion. 
\[
 a_0 = 1,\ a_1 = 3,\ a_2 = 4,\  a_m = a_{m-1} \quad (m \ge 3)
\]
Therefore, the following holds. 
\[
G_e (t) = 1 + 3t + 4t^2 + 4t^3 + 4t^4 + 4t^5 + \cdots
\]

On the other hand, the Cayley graph Cay$(J_3,S)$ is illustrated as follows. 

{
\centering
\begin{tikzpicture}[scale=1.2, >=stealth]
    \tikzset{
        triangle/.style={
            draw, fill=white, line width=0.6pt,
            path picture={
                \draw[fill=white] (90:0.2cm) -- (210:0.2cm) -- (330:0.2cm) -- cycle;
            }
        },
        square/.style={
            draw, fill=white, line width=0.6pt,
            minimum size=0.32cm, inner sep=0pt
        },
        circle_node/.style={
            draw, fill=white, line width=0.6pt,
            circle, minimum size=0.32cm, inner sep=0pt
        }
    }

    \draw[thick] (-0.2, 1) -- (7.4, 1);
    \draw[thick] (-0.2, 0) -- (7.4, 0);

    \draw[dashed, thick] (-0.9, 0.5) -- (-0.3, 0.5);
    \draw[dashed, thick] (7.5, 0.5) -- (8.1, 0.5);

    \foreach \x in {0.8, 1.6, 2.4, 3.2, 4.0, 4.8, 5.6, 6.4} {
        \draw[thick] (\x, 1) -- (\x, 0);
    }


    \node[square] at (0.8, 0) {};

    \node[square] at (1.6, 1) {};
    \draw[fill=white] (1.6, 0.173) -- (1.45, -0.086) -- (1.75, -0.086) -- cycle;

    \draw[fill=white] (2.4, 1.173) -- (2.25, 0.914) -- (2.55, 0.914) -- cycle;
    \node[circle_node] at (2.4, 0) {};

    \node[circle_node] at (3.2, 1) {};

    \draw[fill=white] (4.0, 1.173) -- (3.85, 0.914) -- (4.15, 0.914) -- cycle;
    \node[circle_node] at (4.0, 0) {};

    \node[square] at (4.8, 1) {};
    \draw[fill=white] (4.8, 0.173) -- (4.65, -0.086) -- (4.95, -0.086) -- cycle;

    \node[square] at (5.6, 0) {};

    \foreach \x in {0.8, 1.6, 2.4, 3.2, 4.0, 4.8, 5.6, 6.4} {
        \filldraw (\x, 1) circle (1.5pt);
        \filldraw (\x, 0) circle (1.5pt);
    }

    \node[above=6pt] at (3.2, 1) {$s_{1,3}$};
    \node[below=6pt] at (2.4, 0) {$s_{2,3}$};
    \node[below=6pt] at (3.2, 0) {$e$};
    \node[below=6pt] at (4.0, 0) {$s_{1,2}$};

\end{tikzpicture}
}

\noindent
The elements are represented diagrammatically by their word length. The elements of length 1 are marked with circles ($\circ$), length 2 with triangles ($\triangle$), and length 3 with squares ($\square$). 
The symbol $e$ represents the identity element. 
Clearly, we can obtain the same recurrence relation of the coefficient of the growth series: 
\[
 a_0 = 1,\ a_1 = 3,\ a_2 = 4,\  a_m = a_{m-1} \quad (m \ge 3). 
\]

In the following, we list the cases for degrees $4$, $5$, and $6$. Since we can enumerate elements up to a certain word length by using the normal form for the cactus groups (see \cite{hama2025purecactusgroupsconfiguration}), we can also verify that this enumeration via exhaustive search coincides with the sequence of coefficients of the growth series obtained by Theorem~\ref{growth_of_cactus_gr}. 

\begin{itemize}
    \item $n=4$
    \begin{itemize}
        \item $D_{5,i} \ (i=-1,0, \dots, 3) :0, 1, 5, 5, 0$
        \item $d_{4,k} \ (k=0, \dots, 3) : 1, 6, 10, 5$
        \item power series expansion : 
        \[\displaystyle \frac{(1+t)^3}{t^2 - tx + 1} = 1 + 6t + 20t^2 + 55t^3 + 145t^4 + 380t^5 + \cdots\]
        \item recurrence formula : $a_m = 3a_{m-1} - a_{m-2}$
    \end{itemize}

    \vspace{1em}
    \item $n=5$
    \begin{itemize}
        \item $D_{6,i} \ (i=-1, 0, \dots, 4) :0, 1, 9, 21, 14, 0$
        \item $d_{5,k} \ (k=0, \dots, 4) : 1, 10, 30, 35, 14$
        \item power series expansion : 
        \[\displaystyle \frac{(1+t)^4}{-t^3 + 6t^2 - 6t + 1} = 1 + 10t + 60t^2 + 305t^3 + 1481t^4 + 7116t^5 + \cdots\]
        \item recurrence formula : $a_m = 6a_{m-1} - 6a_{m-2} + a_{m-3}$
    \end{itemize}

    \vspace{1em}
    \item $n=6$
    \begin{itemize}
        \item $D_{7,i} \ (i=-1, 0, \dots, 5) : 0, 1, 14, 56, 84, 42, 0$
        \item $d_{6,k} \ (k=0, \dots, 5) : 1, 15, 70, 140, 126, 42$
        \item power series expansion :
        \[\displaystyle \frac{(1+t)^5}{t^4 - 10t^3 + 20t^2 - 10t + 1} = 1 + 15t + 140t^2 + 1120t^3 + 8554t^4 + 64526t^5 + \cdots\]
        \item recurrence formula : $a_m = 10a_{m-1} - 20a_{m-2} + 10a_{m-3} - a_{m-4}$
    \end{itemize}
\end{itemize}

\section*{Acknowledgments and AI use
}
The author would like to thank his supervisor, Professor Kazuhiro Ichihara. He would also like to thank Carl-Fredrik Nyberg Brodda for the inspiration of this theme. Generative AI tools (Claude and Gemini Pro) were used to assist in searching for and understanding relevant literature, as well as performing auxiliary computations for specific examples.

\bibliographystyle{amsplain}
\bibliography{ref}

\end{document}